\documentclass[12pt,leqno]{amsart}
\usepackage[dvips]{graphicx}
\usepackage{amsfonts}
\usepackage{amssymb}
\usepackage{amsmath}
\usepackage{amscd}
\usepackage{graphicx}
\usepackage{enumitem} 

\def\DATE{\today}
\newtheorem{theorem}{Theorem}
\newtheorem{definition}[theorem]{Definition}

\newtheorem{proposition}[theorem]{Proposition}

\newcommand\C{\mathbb{C}}
\newcommand\R{\mathbb{R}}
\newcommand\g{\mathfrak{g}}

\newcommand\n{\mathfrak{n}}
\newcommand\h{\mathfrak{h}}
\newcommand\K{\mathbb{K}}
\newcommand\Z{\mathbb{Z}}

\newcommand{\mm}{\mathfrak m }

\newcommand\pf{\noindent{\it Proof. }}
\title{Rigid Lie algebras  and algebraicity}
\author[Elisabeth Remm]{ Elisabeth Remm}
\address{E.R: IRIMAS,
        Universit\'e de Haute Alsace, Facult\'e des Sciences et
        Techniques, 4, rue des Fr\`eres Lumi\`ere,
        68093~Mulhouse~cedex, France.}
\email{Elisabeth.Remm@uha.fr}

\begin{document}

\maketitle

\begin{abstract} A finite dimensional complex Lie algebra $\g$ is called rigid if any sufficiently close Lie algebra is isomorphic to it. We prove that this implies that the Lie algebra of inner derivations of $\g$ is an algebraic Lie subalgebra of $gl(\g)$. We show that in general $\g$ is not algebraic. 

\end{abstract}

{\bf Keywords} Lie algebras. Deformations. Rigidity. Algebraic Lie algebras.

\medskip

{\bf MSC:}

\section{Introduction}
The notion of rigidity of Lie algebra  is linked to the following problem:  when does a Lie brackets $\mu$ on a vector space $\g$ satisfy that  every Lie bracket $\mu_1$ sufficiently close to $\mu$ is of the form $\mu_1 = P\cdot \mu $ for some $P \in GL(\g)$ close to the identity? 
A Lie algebra which satisfies the above condition  will be called rigid. The most famous example is the Lie algebra $sl(2,\C)$ of square matrices of order $2$ with vanishing trace. This Lie algebra is rigid, that is any close deformation is isomorphic to it. Let us note that, for this Lie algebra,  there exists a quantification of its universal algebra. This  led to the definition of the famous quantum group $SL(2)$. Another interest of studying the rigid Lie algebras is the fact that there exists, for a given dimension, only a finite number of isomorphic classes of rigid Lie algebras. So we are tempted to establish a classification. This problem has been solved up to the dimension $8$. To continue in this direction, properties must be established on the structure of these algebras. One of the first results establishes an algebricity criterion \cite{Carles}. However, the notion of algebricity which is used is not the classical notion and it includes non-algebraic Lie algebras in the usual sense. The aim of this work is to show that a the Lie algebra is rigid, then its algebra of inner derivations is algebraic.


\section{Rigidity of multiplications on a finite dimensional vector space }

\subsection{Rigidity on the linear space of skew symmetric bilinear map}
Let $E$ be a $n$-dimensional $\K$-vector space, where $\K$ is an algebraically closed field of characteristic 0. We fix a basis $\mathcal{B}=\{ e_1, \cdots , e_n \}$ of $E$. Let $\mu$ be a skew-symmetric bilinear map on $E$ that is:
$$\mu:E \times E \rightarrow E$$
such that  $\mu$ is bilinear and satisfies $\mu(X,Y)=-\mu(Y,X)$ for any $X,Y \in E$ or equivalently $\mu(X,X)=0$ for any $X \in E$ (the field $\K$ is considered of characteristic 0). The structure constants of $\mu$ related to the basis $\mathcal{B}$ are the scalars $X_{ij}^k$, $i,j,k \in [[1,n]]$,  given by 
\begin{eqnarray}
\label{constantes structure}
\mu(e_i,e_j)=\sum_{k=1}^n X_{ij}^k e_k 
\end{eqnarray}
and satisfying
$$X_{ij}^k=-X_{ji}^k, \qquad i,j \in [[1,  n ]].$$
We denote $\mathcal{V}_n$ the $N$-dimensional $\K$-vector space ($N=\frac{n^2(n-1)}{2}$) whose elements are the $N$-uples  $\{  X_{ij}^k, 1 \leq i < j \leq n, k= 1 , \cdots , n  \}.$ It can be considered as an affine space of dimension
$N$. If $\mathcal{SB}il(E)$ is the set of  skew-symmetric bilinear maps
$$\mu: E \times E \rightarrow E,$$ then (\ref{constantes structure}) shows that this vector space is isomorphic to $\mathcal{V}_n$ and we identify these two vector spaces. So, we shall often write $\mu$ for a point of $\mathcal{V}_n$.
 
 Let $GL(E)$ be the algebraic group of linear isomorphisms of $E$. We have a natural action of $GL(E)$ on $\mathcal{SB}il
(E) $ namely
$$\begin{array}{ccc}
GL(E) \times \mathcal{SB}il
(E)  & \rightarrow & \mathcal{SB}il(E) \\
(f,\mu) & \mapsto & \mu_f 
\end{array}$$
where $\mu_f(X,Y)=f^{-1}(\mu(f(X),f(Y))$ for all $X,Y \in E$. This action is translated in an action of $GL(n,\K)$ on $\mathcal{V}_n$ namely 
$$(f=(a_{ij}), (X_{ij}^k))\rightarrow (Y_{ij}^k)$$
 with 
 \begin{eqnarray}
 \label{equation 2}
 \sum_{k=1}^n Y_{ij}^ka_{sk}=\sum_{l,r=1}^n a_{li}a_{rj}X_{lr}^s.
\end{eqnarray}
Since $f=(a_{ij})\in GL(n,\K),$ the vector $(Y_{ij}^k)$ is completely determined by (\ref{equation 2}).

\begin{definition}
The element $\mu=(X_{ij}^k)$ of $\mathcal{V}_n$ is called rigid if its orbit $\mathcal{O}(\mu) = \{ \mu_f, f \in GL(n,\K ) \}$ associated to the action of $GL(n,\K)$ is open in the affine space $\mathcal{V}_n.$
\end{definition}

This notion is a topological notion. Then the considered topology on the affine space is the Zariski-topology. Remember that in this case, the topology coincides with the finer metric topology. Then $\mu$ is rigid if any neighborhood of $\mu$ in $\mathcal{V}_n$ is contained in $\mathcal{O}(\mu).$

\noindent{\bf Example.} We consider $n=2$ and the bilinear map 
$$\mu(e_1, e_2)=X_{12}^1 e_1+ X_{12}^2 e_2$$
with $X_{12}^1$ or $X_{12}^2$ nonzero. It corresponds to a point $\mu=(X_{12}^1,X_{12}^2) \neq (0,0)$ in $\mathcal{V}_2.$ A direct computation gives
$$\mu_f=(Y_{12}^1, Y_{12}^2) \ {\rm with } \ 
\left(
\begin{array}{l}
Y_{12}^1 \\
Y_{12}^2
\end{array}
\right)=\Delta(f) f^{-1}\left(
\begin{array}{l}
X_{12}^1 \\
X_{12}^2
\end{array}
\right)$$
where $\Delta(f)=\det(f).$ 
Then $\mathcal{O}(\mu)=\mathcal{V}_2 \setminus \{ (0,0) \}$  is open in $\mathcal{V}_2$  and the point $\mu$ is rigid.

\medskip

Let $\mu$ be in $\mathcal{V}_n$ and $G_\mu$  be the isotropy subgroup of $\mu$,
$$G_\mu= \{  f \in GL(n,\C) \,  / \, \mu_f=\mu \}.$$

It is a closed subgroup of $GL(n,\K)$ and the orbit $\mathcal{O}(\mu)$ is isomorphic to the homogeneous algebraic space 
$$\mathcal{O}(\mu)=\frac{GL(n,\K)}{G_\mu}.$$
In particular $\mathcal{O}(\mu)$ can be provided with a differentiable manifold structure and $\dim \mathcal{O}(\mu)={ \rm codim} \, G_\mu.$ As a consequence, $\mu \in \mathcal{V}_n$ is rigid if and only if $\dim GL(n,\K)-\dim G_\mu=\dim \mathcal{V}_n=\frac{n^2(n-1)}{2}.$ This implies $\dim G_\mu= \frac{n^2(3-n)}{2}$ and then $n \leq 3$. For $n=2$ the point $\mu=(X_{12}^1,X_{12}^2)=(1,0)$ is rigid. 

\begin{proposition}
If $n\geq 3$, no element $\mu \in \mathcal{V}_n$ is rigid.
\end{proposition}

\pf If $n >3$ then for any $\mu \in \mathcal{V}_n, \ \dim \mathcal{O}(\mu) < \dim \mathcal{V}_n$  and $\mathcal{O}(\mu) $ is not rigid.
If $n=3$, 
from \cite{Remm dim3} $\dim G_\mu \geq 1$ and the orbit of $\mu$ is not rigid. 

Remark that for $n\geq 4,$ a $\mu $ in $ \mathcal{V}_n$ with $\dim G_\mu=0$ can be found (see  \cite{Remm dim3}). But for such an $\mu$, $\dim \mathcal{O}(\mu)< \dim  \mathcal{V}_n$ so $ \mathcal{O}(\mu)$ is not an open set in $\mathcal{V}_n$.

\subsection{Rigidity in stable subsets of $ \mathcal{V}_n$ }

Let $\mathcal{W}$ be an algebraic subvariety of $ \mathcal{V}_n.$ It is defined by a finite polynomial system on $ \mathcal{V}_n.$ We assume that $ \mathcal{W}$   is stable by the action of $GL(n,\mathbb{K})$ on $ \mathcal{V}_n$, that is,
$$\forall \mu \in \mathcal{W}, \ \mu_f \in  \mathcal{W}.$$
It is the case for example for
\begin{itemize}
  \item $ \mathcal{W}= \mathcal{L}_n$  the set of Lie algebra multiplications, that is,
  $$\mu(\mu(X,Y),Z)+\mu(\mu(Y,Z),X)+\mu(\mu(Z,X),Y)=0$$
  for any $X,Y,Z \in E$, or equivalently 
  $$\sum_{l=1}^n X_{ij}^l X_{lk}^s+X_{jk}^lX_{li}^s+X_{ki}^lX_{lj}^s=0$$
  for any $1 \leq i <j\leq n, 1 \leq k \leq n$ and $1\leq s \leq n.$
  \item $ \mathcal{W}= \mathcal{SA}ss_{n}$  the set of skew-symmetric associative multiplications
  $$\mu(\mu(X,Y),Z)-\mu(X,\mu(Y,Z))=0$$
   for any $X,Y,Z \in E$ or equivalently 
  $$\sum_{l=1}^n (X_{ij}^l X_{lk}^s+X_{il}^sX_{jk}^l )$$
  for any $1 \leq i <j\leq j, 1 \leq k \leq n$ and $1\leq s \leq n.$ This set coincides with the set of multiplications satisfying 
  $$\mu(\mu(X,Y),Z)=0,$$
  that is, the subvariety of $\mathcal{L}_n$ constituted of $2$-step nilpotent Lie algebras.
  \item 
  $\mathcal{W}= \mathcal{N}il_n= \{ \mu \in \mathcal{L}_n \ / \ \mu \ {\rm is \ nilpotent}  \} \ {\rm or } \
     \mathcal{S}ol_n= \{ \mu \in \mathcal{L}_n \ / \ \mu \ {\rm is \ solvable}  \}.$
     Recall that $\mu $ is nilpotent if the linear operators 
     $$ad_\mu X: Y \rightarrow \mu(X,Y) $$
     are nilpotent. It is called $k$-step nilpotent if $(ad_\mu X)^{k}=0$ for any $X$ and if there exists $Y$ such that $(ad_\mu Y)^{k-1}\neq0$ ($k$ is also called the nilindex of $\mu $).
 \end{itemize}

\begin{definition}
Let $\mathcal{W}$ be a stable algebraic subvariety  of $\mathcal{V}_n$. An element $\mu \in \mathcal{W}$ is 
$\mathcal{W}$-rigid if the orbit $\mathcal{O}(\mu)$ is open in $\mathcal{W}$.
\end{definition}
Since $\mathcal{W}$ is stable, $\mathcal{O}(\mu) \subset \mathcal{W}$ and 
$$\dim \mathcal{O}(\mu)=n^2-\dim G_\mu.$$
 
\medskip

\noindent{\bf Remark.}  It is also interesting to consider  stable but not necessarily closed subsets $\mathcal{\widetilde{W}}$ of $\mathcal{V}_n$ or of a stable subvariety $\mathcal{W}$ of $\mathcal{V}_n$ that is   for every $\mu \in \mathcal{\widetilde{W}}$ then $\mathcal{O}(\mu) \subset \mathcal{\widetilde{W}}$. For example the subset $\mathcal{N}il_{n,k}$ of $\mathcal{N}il_n$ whose elements are $k$-step nilpotent,  is stable for the action of $GL(n,\K)$.  For this stable subset, there exists another invariant up an isomorphism which permits to describe it: the characteristic sequence  of  a  nilpotent Lie algebras multiplication (see for example \cite{RBreadth} for a detailled presentation of this notion).  Let be $\mu \in \mathcal{N}il_n$. For any $X \in E,$ let $c(X)$ be the ordered sequence, for the lexicographic order, of the dimensions of the Jordan blocks of the nilpotent operator $ad_\mu X$. The characteristic sequence of $\mu$ is the invariant, up to isomorphism,
$$c(\mu)=\max \{c(X), \ X \in E\}.$$ 
In particular, if $c(\mu)=(c_1,c_2,\cdots,1)$, then $\mu$ is $c_1$-step nilpotent.  A vector $X \in E$ such that $c(X)=c(\mu)$ is called a characteristic vector of $\mu$. For a given sequence $(c_1,c_2,\cdots,1)$ with $c_1 \geq c_2\geq \cdots \geq c_s \geq 1$ with $c_1+c_2+\cdots +c_s+1=n$, the subset $\mathcal{N}il_n^{(c_1,c_2,\cdots,1)}$ is stable in  $\mathcal{N}il_n$.

\begin{definition}
Let $\mathcal{\widetilde{W}}$ be a stable subset of $\mathcal{V}_n$ or of a stable subvariety $\mathcal{W}$ of $\mathcal{V}_n$. A multiplication $\mu \in \mathcal{\widetilde{W}}$ is called $\mathcal{\widetilde{W}}$-rigid if for any neighborhood $V(\mu)$ of $\mu$ in  $\mathcal{V}_n$ or  in $\mathcal{W}$ then $V(\mu) \cap \mathcal{V}_n$ or $V(\mu) \cap \mathcal{W} \subset \mathcal{O}(\mu)$.
\end{definition}

This notion of $\mathcal{\widetilde{W}}$-rigidy has been introduced in \cite{GR-Kegel} to study  the set of $k$-step nilpotent Lie algebras.
 \subsection{How to prove the rigidity}
 We have two approaches 
 
 \noindent 1. {\bf A topological way}.  By definition, $\mu $ is rigid if $\mathcal{O}(\mu)$ is open in $ \mathcal{W} $ and $\mathcal{W} \setminus \mathcal{O}(\mu)$ is an algebraic subset of $\mathcal{W}.$ But $\mathcal{O}(\mu)$ is provided with a differentiable homogeneous manifold contained in the affine space $\mathcal{V}_n.$ So we can consider open neighbourhood of $\mu$ for the "metric" topology. Now, although $ \mathcal{W}$ contains singular points, any rigid point $\mu$ in $ \mathcal{W}$ is non singular since its orbit is open in $ \mathcal{W}$. Thus in order to prove the rigidity we can consider an open neighbourhood $B_\mu$ of $\mu$ in  $\mathcal{V}_n$. For this we can use a method inspired by the determination of the algebraic Lie algebra of an algebraic Lie group using the dual numbers and consider non archimedian extension of $\K$.

 \noindent 2. {\bf A geometrical way}. Since $ \mathcal{O}(\mu)$ is a differentiable manifold its tangent space $T_\mu  \mathcal{O}(\mu)$ to $\mu$ is well defined. It is isomorphic to the quotient space $\frac{gl(n,\K)}{\mathcal{D}er(\mu)}$ where $\mathcal{D}er(\mu)=\{ f \in gl(n,\K), \delta f(X,Y)=\mu(f(X),Y)+\mu (X,f(Y))-f(\mu(X,Y))=0 \}$ is the algebraic Lie algebra of $G_\mu.$  We deduce 
 $$\dim T_\mu \mathcal{O}(\mu)=n^2-\dim \mathcal{D}er (\mu)$$
 and $T_\mu \mathcal{O}(\mu)$ is isomorphic to the subspace of bilinear maps whose elements are $\delta f$ for any $f \in gl(n,\C)$ generally denoted $B^2(\mu, \mu).$ Then
\begin{proposition}
The application  $\mu$ is rigid in $\mathcal{W}$ if
 $$\dim T_\mu \mathcal{O}(\mu)=\dim B^2(\mu, \mu) =\dim T_\mu \mathcal{W}.$$
\end{proposition}

 The determination of $T_\mu \mathcal{W}$ is a little bit difficult. Since we assume that $\mu$ is rigid, necessarily $T_\mu \mathcal{W}$ exists. In a first time we can compute the Zariski tangent space $T^Z_{\mu} \mathcal{W}$  usually denoted $Z^2(\mu, \mu),$ defined by the linear system obtained by considering the polynomial system of definition of $\mathcal{W}$, translated to the point $\mu$ and taking its linear part. But $T_\mu \mathcal{W}\subset T^Z_\mu \mathcal{W}=Z^2(\mu, \mu)$ and this inclusion can be strict. This appears as soon as the affine schema  which defines $\mathcal{W}$ is not reduced at the point $\mu$. We can illustrate this in a simple  example. Let us consider the algebraic variety $M $ in $\C^3$ defined by the polynomial system
 $$\left\{ 
 \begin{array}{l}
 X_1 X_2 + X_1 X_3-2X_2X_3=0, \\
 2X_1 X_2 -3 X_2 X_3+ X_3^2=0 
 \end{array}
 \right.
$$
Only the point $(0,0,0)$ is singular. Let us compute $T_\mu M$ at the point $(1,1,1)$. Linearizing the system we obtain 
$$(\star) \left\{ 
 \begin{array}{l}
 2X_1 -X_2 - X_3+X_1 X_2 +X_1 X_3 -2X_2X_3=0, \\
 2X_1 -X_2 - X_3+2X_1 X_2 -3X_2 X_3+ X_3^2=0 
 \end{array}
 \right.
$$
and $T_\mu M=Ker\{ \rho(X_1,X_2,X_3)=2X_1-X_2-X_3 \}$ ant it is of dimension 2 although $M$ is a one-dimensional curve. To compute $T_\mu M$ we come back to the definition of the tangent vector. We consider a point 
$(\varepsilon_1, \varepsilon_2, \varepsilon_3)$ in $M$ close to $(0,0,0)$ and satisfying $(\star).$ Thus 
$$\left\{ 
 \begin{array}{l}
 2\varepsilon_1 -\varepsilon_2 - \varepsilon_3+\varepsilon_1 \varepsilon_2 +\varepsilon_1 \varepsilon_3 -2\varepsilon_2\varepsilon_3=0, \\
 2\varepsilon_1 -\varepsilon_2 - \varepsilon_3+2\varepsilon_1 \varepsilon_2 -3\varepsilon_2 \varepsilon_3+ \varepsilon_3^2=0 
 \end{array}
 \right.
$$
A similar approach with the dual number implies
$$2\varepsilon_1-\varepsilon_2-\varepsilon_3=0$$
that is $(\varepsilon_1,\varepsilon_2,\varepsilon_3)\in T_\mu^2 M$ and 
$$\left\{ 
 \begin{array}{l}
\varepsilon_1 \varepsilon_2 +\varepsilon_1 \varepsilon_3 -2\varepsilon_2\varepsilon_3=0, \\
 2\varepsilon_1 \varepsilon_2 -3\varepsilon_2 \varepsilon_3+ \varepsilon_3^2=0 
 \end{array}
 \right.
$$
This is equivalent to 
$$\varepsilon_1= \frac{2\varepsilon_2\varepsilon_3}{\varepsilon_2+\varepsilon_3}=\frac{\varepsilon_3(3\varepsilon_2-\varepsilon_3)}{2\varepsilon_2}$$ that is $$\varepsilon_2=\varepsilon_3=\epsilon_1$$ and $T_\mu W = \{ (x, x,x), x \in \C \} \varsubsetneqq T_\mu ^Z W.$

Remarks. 1) We know many examples of rigid Lie algebras such that $T_\mu \mathcal{O}(\mu) \varsubsetneqq T_\mu ^Z \mathcal{W}.$ Since $T_\mu \mathcal{O}(\mu)$ and  $T_\mu ^Z \mathcal{W}$  coincide respectively with the space of coboundaries and the space of $2$-cocycles associated with the Chevalley Eilenberg cohomology $H^*(\mu, \mu)$ of $\mu$ ( when $\mu$ is a Lie algebra) the classical theorem $\dim H^2(\mu, \mu)=0$ implies $\mu$ is rigid. But the converse is not true because the determination of $T_\mu ^Z \mathcal{W}$ is not suffisant to compute $\dim T_\mu \mathcal{W}.$ There exists another approach of the rigidity using the deformation theory close to the cohomogical point of view. We consider a formal series $\mu_t=\mu + t \varphi_1 + \cdots + t^n \varphi_n +\cdots $ and $\mu_t \in \mathcal{W}$ implies $\varphi_1 \in T_\mu ^Z W.$ But to compute $T_\mu  \mathcal{W}$ it is necessary to look all the relations between $\varphi_1, \varphi_2, \cdots$

2) As the elements of $T_\mu ^Z \mathcal{W}$ can be interpreted  as cocycle associated with the Chevalley Eilenberg cohomology, the elements of $T_\mu \mathcal{W}$ can be interpreted as particular cocycles. This will be the aim of a next work.

\subsection{Consequence of the reductivity of $GL(n,\K)$}

The action of $GL(n,\K)$ on $\mathcal{V}_n$ or on an algebraic subvariety $\mathcal{W}$  is an example of an action of reductive group on an algebraic affine variety. Recall that any element $f \in GL(n,\K)$ decomposes as $f=f_s \circ f_u$ where $f_s$ is a semi-simple and $f_u$ is unipotent. So the fact that $G$ is reductive implies that the maximal normal unipotent subgroup $R_u(G)$ of $G$ is trivial. Let $G_U$ be the maximal unipotent subgroup of $G$. 

\begin{proposition}
Let $\mu $ be in $\mathcal{L}ie_n$ and $G_\mu$ the maximal unipotent subgroup of $G$. Then the orbit $\mathcal{O}_{G_\mu}(\mu)=\{ \mu_f \, / \, f\in G_\mu\}$ is closed in  $\mathcal{L}ie_n$.
\end{proposition}

\noindent{\bf Consequence.} As we are interested in open orbits, we will be concerned with the action of the maximal torus (all the elements are semi-simple and commuting). Let $\mu=(X_{ij}^k ) $ be in  $\mathcal{L}ie_n$ and $f\in T$ where $T$ is a maximal torus. We can suppose that the basis $\{ X_1, \cdots , X_n\}$ associated to the  $X_{ij}^k$'s is a basis of eigenvectors of $f$. So 
$$\mu_f(X_i,X_j)=\sum_k
 \frac{\lambda_i \lambda_j}{\lambda_k}X_{ij}^k
X_k$$
and the structure constants $Y_{ij}^k$ of $\mu_f$ satisfy 
$$Y_{ij}^k= \frac{\lambda_i \lambda_j}{\lambda_k}X_{ij}^k.$$
We then consider a point $\tilde{\mu}$ in a neighborhood  of $U$ defined by 
$$ Z_{ij}^k=X_{ij}^k(1+\rho_{ij}^k).$$
The rigidity of $\mu$ implies that for $\rho_{ij}^k \simeq 0$, there exist $\{ \lambda_i \}$ such that 
$$ \rho_{ij}^k=\frac{\lambda_i \lambda_j-\lambda_k}{\lambda_k}.$$
We will come back to this system in the Rank Theorem \cite{A.G-rank}.

\section{Rigid Lie algebras and algebraicity}

\subsection{Rigid Lie algebras}
Recall that a $\K$-Lie algebra is a pair $\g=(V,\mu)$ where $\g$ is a $\K$-vector space and $\mu$ a Lie algebra multiplication on $V$. In this section we need to differentiate $\g$ to $\mu$ because the structure of $\g$ depends to the nature of the vector space $V$.

\begin{definition}
A finite dimensional $\K$-Lie algebra $\g=(V,\mu)$ is called rigid if $\mu$ is $\mathcal{L}ie_n$-rigid, where $n=\dim V$.
\end{definition}

From the previous discussion, a Lie algebra is rigid if and only if $\dim T_\mu \mathcal {O}(\mu)= \dim T_\mu \mathcal{L}ie_n$. In particular, since the Zariski tangent space $T^Z_\mu \mathcal{L}ie_n$ contains $T_\mu \mathcal{L}ie_n$, we have the classical Nijenhuis-Richardson theorem that we can write:

\begin{theorem}({\rm Nijenhuis-Richardson}) 
A $\K$-Lie algebra $\g=(V,\mu)$ satisfying $\dim H^2(\g,\g)=0$ is rigid.
\end{theorem}
In fact $H^2(\g,\g)= \frac{Z^2(\g,\g)}{B^2(\g,\g)}$ and $Z^2(\g,\g)$ is isomorphic to $T^Z_\mu \mathcal{L}ie_n$ and $B^2(\g,\g)$ is isomorphic to $T_\mu \mathcal{O}(\mu)$.

\noindent {\bf Remarks.} 

1. We have discussed the converse in the previous section. The simpliest example of rigid Lie algebra having a non zero-dimensional $H^2(\g,\g)$  actually known is in dimension $13$. It is given by the multiplication
$$
\left\{
\begin{array}{l}
\lbrack X_1,X_i \rbrack=X_{i+1}, 2\leq i \leq 11, \  \lbrack X_2,X_i \rbrack=X_{i+2}, \ 3\leq i \leq 10, \\
\lbrack T ,X_i \rbrack=iX_i, 1\leq i \leq 12. 
\end{array}
\right.
$$
Here the dimension of the second space of cohomology is $1$. This means that the $$\dim T^Z_\mu(\mathcal{L}ie_{13})-\dim T_\mu (\mathcal{L}ie_{13})=1.$$ A direct computation considering a point of $\mathcal{L}ie_{13}$ close to $\mu$ shows that 
$$T^Z_\mu(\mathcal{L}ie_{13})=T_\mu(\mathcal{L}ie_{13}) \oplus \C{\phi}$$
where $\phi$ is the element of $\mathcal{V}_{13}$ given by
$$
\left\{
\begin{array}{l}
\phi(X_2,X_i)=(4-i)X_{2+i}, \ 5 \leq i \leq n-2,\\
\phi (X_3,X_i)=X_{3+i}, \ 4 \leq i \leq n-3.\\
\end{array}
\right.
$$
This non tangent "cocycle" has been already defined in \cite{Nico}.

2.  Assume that $\mu$ is rigid in $\mathcal{L}ie_n$. In this case $T_\mu(\mathcal{L}ie_n)=T_{\mu_1}(\mathcal{L}ie_n)$ for any $\mu_1 \in \overline{\mathcal{O}(\mu)}$. For example, if the rigid Lie algebra $\g=(V,\mu)$ is a contact rigid Lie algebra, then $T_\mu(\mathcal{L}ie_n)$ can be computed considering the tangent space at $\mu_1$ where $\mu_1$ is the multiplication of the Heisenberg algebra.

\subsection{Algebraic Lie algebras}
Recall that $\K$ is an algebraically closed field of characteristic $0$. (A study of rigid Lie algebra when $\K=\R$ is proposed in \cite{ACG}).

\begin{definition} \label{aa}
A $\K$-Lie algebra is algebraic if it is the Lie algebra of an algebraic Lie group. 
\end{definition}

\noindent Examples: \begin{enumerate}
  \item Any complex semi-simple Lie algebra is algebraic.
  \item Any Lie algebra which coincides with its derived sub-algebra is algebraic.
  \item For any Lie algebra $\g$, there exists an algebraic Lie algebra, containing $\g$ and having the same derived subalgebra. It is called the algebraic Lie algebra generated by $\g$.
\end{enumerate}

\noindent Problem: A Lie algebra is not always algebraic so how to characterize an algebraic Lie algebra? 
There exists some criterium to study the algebraicity of a Lie algebra. We can always assume that an algebraic Lie algebra is linear that is it is a Lie subalgebra of some $gl(n,\C)$.  Let  $\g$ be a Lie sub-algebra of $gl(n,\K)$. A replica $Y$ of an element $X \in \g$ is an element of the algebraic sub-algebra $\g(X)$ which is the smaller algebraic Lie sub-algebra of $gl(n,\K)$ containing $X$. A Lie sub-algebra $\g$ of $gl(n,\K)$ is algebraic if and only if for any $X \in \g$, all the replica of $X$ are in $\g$.

Recall also the structure of algebraic solvable or nilpotent Lie algebras.

\begin{proposition}
Let  $\g$ be an algebraic nilpotent Lie algebra, subalgebra of $gl(n,\K)$. Let  $\frak{n}$ be the ideal of $\g$ whose elements are the nilpotent elements of $\g$. Then $\g$ is the direct sum $$\g=\frak{n} \oplus \frak{a}$$
where $\frak{a}$ is an abelian algebraic subalgebra of $\g$ contained in the center of $\g$ whose all the elements are semi-simple.
\end{proposition}
Let us note that  $\frak{a} \subset Z(\g)$. In fact if  $X \in \frak{a}$, then $X$ semi-simple implies that $ad X$ is also semi-simple. But it is also nilpotent because  $\g$ is nilpotent, then  $ad X=0$ et $X \in Z(\g)$. 

Concerning the solvable case, we have
\begin{proposition}
Let  $\g$ be an algebraic solvable Lie algebra, subalgebra of $gl(n,\K)$. Let  $\frak{n}$ be the ideal of $\g$ whose elements are the nilpotent elements of $\g$. Then $\g$ is the direct sum $$\g=\frak{n} \oplus \frak{a}$$
where $\frak{a}$ is an abelian algebraic sub-algebra of $\g$ with only  semi-simple elements.
\end{proposition}

\medskip

\noindent{Examples.}

\begin{enumerate}
\item 
The one-dimensional abelian Lie algebra
$$\frak{a_1}= 
\left\{\begin{pmatrix}
   x   & x  \\
    0  &  x\\
\end{pmatrix}, \ x \in \C \right\}.$$
is not algebraic. In fact the semisimple part $ 
\left\{\begin{pmatrix}
   x   & x  \\
    0  &  x\\
\end{pmatrix}
\right.$ is not in $\frak{a_1}$ as soon as $x \neq 0.$

\item Let us consider the following $3$-dimensional Lie algebras:
 $$\frak{n}_1=\left\{ 
\left(
\begin{array}{llll}
x_1+x_2 & x_1+x_2 & 0 & x_1 \\
x_1+x_2 & x_1+x_2 & 0 & x_2 \\
 x_1& x_1+2 x_2 & 0 & x_3 \\
0 & 0 & 0 & 0 
\end{array}
\right), \ x_i \in \C \right\} \ \ 
\frak{n_2}= 
\left\{\begin{pmatrix}
   0   & x_1 & x_3   \\
    0  &  0 & x_2 \\
    0 & 0 & 0
\end{pmatrix}, \ x_i \in \C \right\}.$$
They are isomorphic as Lie algebras but $\n_2$ is algebraic and $\n_1$ not (see \cite{RGarxiv}). We can even construct a family of $3$-dimensional non algebraic Lie algebras which are isomorphic, as Lie algebra, to the algebraic Lie algebra $\n_2$. We consider the $3$-dimensional linear subspace $\h_{\alpha,\beta}$ of $gl(4,\C)$ 
$$\left\{ 
\left(
\begin{array}{llll}
x_1+x_2 & x_1+x_2 & 0 & x_1 \\
x_1+x_2 & x_1+x_2 & 0 & x_2 \\
\alpha x_1+(\beta -1) x_2 & \beta x_1+(\alpha +1) x_2 & 0 & x_3 \\
0 & 0 & 0 & 0 
\end{array}
\right), \ x_1,x_2,x_3 \in \C \right\}$$
where  $\alpha, \beta$ are given elements of $\C$. A basis is given by
$$X_1=
\left(
\begin{array}{llll}
1 & 1 & 0 & 1 \\
1 & 1 & 0 & 0 \\
\alpha  & \beta  & 0 & 0 \\
0 & 0 & 0 & 0 
\end{array}
\right)
, \ X_2=
\left(
\begin{array}{llll}
1 & 1 & 0 & 0 \\
1 & 1 & 0 & 1 \\
\beta -1  & \alpha +1 & 0 & 0 \\
0 & 0 & 0 & 0 
\end{array}
\right), \ X_3=\left(
\begin{array}{llll}
0 & 0 & 0 & 0 \\
0 & 0 & 0 & 0 \\
0  & 0 & 0 & 1 \\
0 & 0 & 0 & 0 
\end{array}
\right)
$$
We verify that
$$[X_1,X_2]=X_3, [X_1,X_3]=[X_2,X_3]=0.$$
Then  $\h_{\alpha,\beta}$ is a  $3$-dimensional Lie subalgebra of  $gl(4,\C)$. It is isomorphic to the $3$-dimensional Heisenberg Lie algebra, but this isomorphism is not the linear part of isomorphism of algebraic groups. 
Let us consider $X_1$ with $\alpha+\beta \neq 0$. Its Chevalley-Jordan decomposition  
$$X_1=X_{1,s}+X_{1,n}$$
is given by
 $$X_{1,s}=
\left(
\begin{array}{llll}
1 & 1 & 0 & 1/2 \\
1 & 1 & 0 & 1/2 \\
\frac{\alpha+\beta}{2}  &\frac{\alpha+\beta}{2}   & 0 &\frac{\alpha+\beta}{4} \\
0 & 0 & 0 & 0 
\end{array}
\right), \ X_{1,n}=
\left(
\begin{array}{llll}
0 & 0 & 0 & 1/2 \\
0 & 0 & 0 & -1/2 \\
\frac{\alpha-\beta}{2}   & -\frac{\alpha-\beta}{2}  &0 & -\frac{\alpha+\beta}{4} \\
0 & 0 & 0 & 0 
\end{array}
\right)
$$
Since $X_{1,s} \notin \h_{\alpha,\beta}$ and also $X_{1,n} \notin \h_{\alpha,\beta}$, we deduce
\begin{proposition}
The $3$-dimensional nilpotent lie subalgebra  $\h_{\alpha,\beta}$ of $gl(4,\C)$ is not algebraic.
\end{proposition}

A matrix $X \in  \h_{\alpha,\beta}$ is nilpotent if and only if  $x_1+x_2=0$. In fact the eigenvalues of $X$ are $0$ and $2(x_1+x_2)$. We deduce that the set of nilpotent  matrices of  $\h_{\alpha,\beta}$ is the subspace generated by  $X_1-X_2,X_3$ and it is an abelian ideal of dimension $2$. Let us note that any non trivial matrix of  $\h_{\alpha,\beta}$ is diagonalisable and   $\h_{\alpha,\beta}$ doesn't admit a Chevalley decomposition.

\medskip

We have recalled that any Lie algebra $\g_0$ generates an algebraic Lie algebra $\g_1$ which is the smallest algebraic algebra containing $\g_0$ and these two algebras are the same derived Lie algebra. Let us determinate the algebraic Lie algebra generated by  $\h_{\alpha,\beta}$. This algebra contains the semi-simple part for any $X \in  \h_{\alpha,\beta}$. If $X=x_1 X_1+ x_2 X_2+ x_3 X_3$ 
 $$X_{s}=
\left(
\begin{array}{llll}
x_1+x_2 & x_1+x_2 & 0 & \frac{x_1+x_2}{2} \\
x_1+x_2 & x_1+x_2 & 0 & \frac{x_1+x_2}{2} \\
(x_1+x_2 )\frac{\alpha+\beta}{2}  &(x_1+x_2 )\frac{\alpha+\beta}{2}   & 0 &(x_1+x_2 )\frac{\alpha+\beta}{4} \\
0 & 0 & 0 & 0 
\end{array}
\right)$$
Consider $X_4=\left(
\begin{array}{llll}
1 & 1 & 0 & 1/2 \\
1 & 1 & 0 & 1/2 \\
\frac{\alpha+\beta}{2}  &\frac{\alpha+\beta}{2}   & 0 &\frac{\alpha+\beta}{4} \\
0 & 0 & 0 & 0 
\end{array}
\right)$. Then we have 
$[X_i,X_4]=0$ for $i=1,2,3$.
and  $\h_{\alpha,\beta} \oplus \K\{X_4\}$ is a 4-dimensional Lie algebra, containing  $\h_{\alpha,\beta}$. Moreover, for any $X$ in 
$\mm$, $X_s$ and $X_n$ the semi-simple and nilpotent parts of the Jordan decomposition of $X$ are  in $\mm$.
This Lie algebra is 
$$\mm= \left\{ 
\left(
\begin{array}{llll}
x_1+x_2+x_4 & x_1+x_2+x_4 & 0 & x_1+x_4/2 \\
x_1+x_2+x_4 & x_1+x_2+x_4 & 0 & x_2 +x_4/2\\
\alpha x_1+(\beta -1) x_2+\frac{\alpha+\beta}{2} x_4 & \beta x_1+(\alpha +1) x_2+ \frac{\alpha+\beta}{2} x_4 & 0 & x_3+\frac{\alpha+\beta}{4}x_4  \\
0 & 0 & 0 & 0 
\end{array}
\right)
\right\},$$
with $ x_1,x_2,x_3,x_4 \in \C. $
where  $\alpha, \beta$ are given elements of $\C$. Let us note that an element of $\mm$ is nilpotent if and only if $x_1+x_2+x_4=0.$ Then the set $\n_1$ of nilpotent elements of $\mm$ is the $3$-dimensional linear subspace of $\mm$
$$\n_1= \left\{ 
\left(
\begin{array}{llll}
0 & 0 & 0 & x_1/2-x_2/2 \\
0& 0 & 0 & -x_1 /2+x_2/2\\
\frac{\alpha-\beta}{2}  x_1+\frac{-\alpha+\beta-2}{2} x_2 & \frac{-\alpha+\beta}{2} x_1+\frac{\alpha-\beta+1}{2} (\alpha +1) x_2+ & 0 & x_3+\frac{\alpha+\beta}{4}(-x_1-x_2)  \\
0 & 0 & 0 & 0 
\end{array}
\right)
\right\},$$
The set of  diagonalisable elements is the $1$-dimensional susalgebra
$$\frak{a}_1= \left\{ 
\left(
\begin{array}{llll}
y & y & 0 & y/2 \\
y & y & 0 & y/2\\
\frac{\alpha+\beta}{2} y & \frac{\alpha+\beta}{2} y & 0 & \frac{\alpha+\beta}{4}y  \\
0 & 0 & 0 & 0 
\end{array}
\right)
\right\},$$
with $y \in \C$ and  $\mm$ is decomposable, that is
$$\mm= \n_1 \oplus \frak{a}_1.$$
Moreover, for any $X \in \mm$, its components $X_s$ and $X_n$ are also in $\mm$.  To study the algebraicity of $\mm$ we have to compute for any $X \in \mm$, the algebraic Lie algebra $\g(X)$ generated by $X$. Let $X_s$ its semisimple component. The eigenvalues are $0$  which is a  triple root and $ 2(x_1+x_2+x_4 )$. We assume that $x_1+x_2+x_4  \neq 0$. The set $\Lambda$ is constituted of $4$-uples of integers $(p_1,p_2,p_3,0)$. If $Y$ is a semisimple element of $\mm$ commuting with $X_s$, its eigenvalues $(\mu_1,\mu_2,\mu_3,\mu_4$ satisfy $p_1\mu_1+p_2\mu_2+p_3\mu_3=0$ for any $p_1,p_2,p_3 \in \Z$. Then $\mu_1=\mu_2=\mu_3=0$. We deduce that
$$Y=
\left(
\begin{array}{cccc}
 \frac{m}{2} & \frac{m}{2} & 0 & \frac{m}{4} \\
 \frac{m}{2} & \frac{m}{2} & 0 & \frac{m}{4} \\
 \frac{m (\alpha +\beta )}{4 } &  \frac{m (\alpha +\beta )}{4 } & 0 &  \frac{m (\alpha +\beta )}{8 } \\
 0 & 0 & 0 & 0 \\
\end{array}
\right)
$$
and $\g(X_s) \subset \mm$. Let  $X_n=X-X_s$ be the nilpotent component of $X$. Then 
$$X_n=\left(
\begin{array}{llll}
0 & 0 & 0 & \frac{x_1-x_2}{2} \\
0 & 0 & 0 & \frac{-x_1+x_2}{2} \\
\frac{\alpha-\beta}{2} x_1+(\frac{-\alpha+\beta}{2}-1)x_2   & \frac{-\alpha+\beta}{2} x_1+(\frac{\alpha-\beta}{2}+1)x_2  & 0 &x_3-(x_1+x_2 )\frac{\alpha+\beta}{4} \\
0 & 0 & 0 & 0 
\end{array}
\right)
$$
Such vector belongs to an algebraic nilpotent $3$-dimensional Lie algebra whose all its elements are nilpotent. Since $\g(X_n)$ is contained in this algebra, it is also contained in $\mm$. Then we have

\begin{proposition}
The Lie algebra $\mm$ is algebraic. It is the algebraic Lie algebra generated by $\h_{\alpha,\beta}$.
\end{proposition}

\end{enumerate}

\subsection{Rigidity and algebraicity}

Recall that in \cite{Carles} we have the following result

\begin{center} 
{\it Any complex rigid Lie algebra is algebraic.}
\end{center}

In this form this result is not true. In fact,
let us consider the following $2$-dimensional Lie algebras
 $$\g_1=\left\{ \begin{pmatrix}
     x &  y  \\
    0  &  0
\end{pmatrix}, \ x, y \in \C
\right\}, \ \ \g_2=
\left\{
\begin{pmatrix}
  x	&   x &y \\
   0   &  x & 0 \\
   0 & 0 & 0
\end{pmatrix}, \ x, y \in \C
\right\}.
$$
These Lie algebras are isomorphic, they have the same Lie multiplication
$$\mu(X,Y)=Y.$$
Computing the  replica of any elements of $\g_1$, we can conclude that this Lie algebra is algebraic.  Concerning $\g_2$, the semi-simple part of the element  corresponding to $x=1$ and $y=0$ is not in $\g_2$. This implies that $\g_2$ is not algebraic.

\begin{proposition}
These two Lie algebras $\g_1$ and $\g_2$ are rigid. But $\g_2$ is a rigid non algebraic Lie algebra.
\end{proposition}

The mistake in the Carles's result is not a consequence of a bad computation but to a strange definition of the algebraicity. In its paper, a Lie algebra is called algebraic if it is isomorphic as a Lie algebra to a Lie algebra of an algebraic group. If we consider the previous counter example, the Lie algebra $\g_2$ is isomorphic to the algebraic Lie algebra $\g_1$ and it is algebraic in the Carles's sense, but not algebraic in the classical sense. The isomorphism between $\g_2$ and $\g_1$ is a Lie algebras isomorphism but not an algebraic Lie algebras isomorphism. From the Carles's definition, we deduce that
any nilpotent Lie algebra is algebraic. This is wrong in general. We have given examples in the previous section.
\medskip

\noindent Thus we consider in the following the classical definition of the algebraicity. In this classical context, we shall prove that if $\mu$ is a rigid Lie algebra multiplication, then the Lie algebra $\mathcal{A}d_\mu$ whose elements are the linear operators $ad_\mu X$ is algebraic.
 
  \subsection{The Lie algebra $\mathcal{A}d_\mu$}

 We have seen that the notion of rigidity is combined with the Lie algebra multiplication and not with the Lie algebra. The notion of rigidity of a Lie algebra is given in terms of rigidity of its Lie multiplication. But a Lie multiplication on $n$-dimensional vector space $E$ defines a natural subalgebra of $gl(E)$ that is
 the Lie algebra $\mathcal{A}d_\mu$ whose elements are the operators $ad_\mu X$ for any $X \in E$.  Let us recall some classical results. Let $\mathcal{D}er_\mu$ be the Lie algebra of derivations of $\g$. It is an algebraic Lie subalgebra of $gl(E)$. The Lie algebra $\mathcal{A}d_\mu$ is an ideal of  $\mathcal{D}er_\mu$ but it is not in general an algebraic Lie subalgebra of $\mathcal{D}er_\mu$. For example, let us consider the $3$-dimensional Lie multiplication given by
 $$[T,X_1]=eX_1,\ [T,X_2]=\pi X_2.$$
 The element $ad T$ is semi-simple with $0$, $e$ and $\pi$ as eigenvalues. Let $(m_1,m_2,m_3)$ be in $Z^3$  such that $m_1 0+m_2 e+m_3\pi=0$. Then $m_2=m_3=0$ and any replica of $ad T$ is semi-simple with eigenvalues $\lambda_1,\lambda_2,\lambda_3$ satisfying $$m_1\lambda_1+m_2\lambda_2+m_3\lambda_3=m_1\lambda_1=0$$
 that is $\lambda_1=0.$ Then the replica corresponding to $\lambda_1=0,\lambda_2=\lambda_3 \neq 0$ is not in $\mathcal{A}d_\mu$. Then this linear Lie algebra $\mathcal{A}d_\mu$ is not algebraic. 
 
\medskip

 We know some general situation where $\mathcal{A}d_\mu$ is an algebraic linear Lie algebra. For example
 \begin{enumerate}
  \item If $\mu$ is a nilpotent Lie multiplication on $E$, then from Engel's theorem, any $ad X$ is nilpotent and the Lie algebra $\mathcal{A}d_\mu$ is unipotent. This implies that $\mathcal{A}d_\mu$ is algebraic.
  \item If any derivation of $\mu$ is inner, that is $\mathcal{A}d_\mu=\mathcal{D}er_\mu $ (that is the first space of the Chevalley-Eilenberg of $\mu$ is trivial), then $\mu$ is algebraic. 
  \item Let $\g \subset gl(V)$ be an algebraic linear Lie algebra. If $\mu$ is the Lie multiplication of $\g$ corresponding to the bracket in $gl(V)$, then $\mathcal{A}d_\mu$ is also algebraic. To prove the algebraicity of $ \mathcal{A}d_\mu$, we have to show that for any $ad X \in  \mathcal{A}d_\mu$, $X \in \g$, the set $ \mathcal{A}d_\mu(ad X)$ of replica  of $ad X$ is contained in $ \mathcal{A}d_\mu$. Let us determine this set. Assume for instance that $\g \subset gl(V)$ is a linear Lie algebra (algebraic or not). For any $u \in gl(V)$ which satisfies
 $$[u ,\g] \subset \g$$
 (where $[,]$ is the Lie bracket in $gl(V)$ corresponding to $\mu$ in $\g$),
 we consider the endomorphism of $\g$
 $$\rho_u (X)= [u,X].$$
This is a derivation of $\g$. The sub-algebra $\h$ of $gl(V)$ given by
 $$\h=\{ u \in gl(V), \ \rho_u(\g) \subset \g\}$$
 is an algebraic sub-algebra of $gl(V)$ containing $\g$. Let us denote also by $\widetilde{\g}$ the algebraic Lie algebra generated by $\g$. If $\g$ is not algebraic, then $\g$ is strictly contained in $\widetilde{g}$ but these two Lie algebras have the same derived sub-algebra. Since $\h$ is an algebraic Lie algebra containing $\g$, then
 $$\widetilde{g} \subset \h.$$ We have
 $$ \mathcal{A}d_\mu(ad X)=\{\rho_{\widetilde{X}}, \ \widetilde{X} \in \widetilde{g}\}.$$
 If we assume now that $\g$ is algebraic, then $\g=\widetilde{g}$ and
 $$ \mathcal{A}d_\mu(ad X)=\{\rho_{\widetilde{X}}, \ \widetilde{X} \in \widetilde{g}\}=\{\rho_{X}, \ X\in \widetilde{g}\}=\mathcal{A}d_\mu.$$
\end{enumerate}
 
\noindent Remark:  The converse of this proposition  is not true. For example, if we consider the $2$-dimensional solvable algebra $\g_2$ defined above, we have seen that it is a not algebraic Lie algebra. Let us compute for this Lie algebra  $ \mathcal{A}d_\mu$.
If we consider the basis of this Lie algebra given by $X$ corresponding to $x=1,y=0$ and $Y$ corresponding to $x=0,y=1$, we have $\mu(X,Y)=Y$ and $ \mathcal{A}d_\mu$ is the Lie algebra
 $$
\left\{
\begin{pmatrix}
   0   &  0 \\
   -y & x 
\end{pmatrix}, \ x, y \in \C
\right\}.
$$
 It is an algebraic abelian Lie algebra because all the replica of any element of this Lie algebra are inside.

\medskip

 \begin{theorem}
 If the Lie multiplication $\mu$ is rigid in $\mathcal{L}ie_n$, then the Lie algebra  $ \mathcal{A}d_\mu$ is algebraic.
 \end{theorem}
 \pf  A linear Lie algebra is algebraic if and only if its radical is algebraic. Thus we can assume that $\mathcal{A}d_\mu$ is solvable. Let be $X \in E$ and $U=ad X$. Let $\widetilde{\mathcal{A}d_\mu}$ the algebraic Lie algebra generated by $\mathcal{A}d_\mu$. Since $\mathcal{A}d_\mu$ is a Lie sub-algebra of the algebraic Lie algebra $\mathcal{D}er_\mu$, then $\widetilde{\mathcal{A}d_\mu}$ is a Lie algebraic sub-algebra of $\mathcal{D}er_\mu$, and  the replica $\widetilde{U}$ of $U$ belongs to $\mathcal{D}er(\mu)$ and it is a derivation of $\mu$. We can assume that $\mu$ is not a nilpotent Lie multiplication because as we have recalled above $ \mathcal{A}d_\mu$ is algebraic. This implies that the derivation $\widetilde{U}$ is singular. Moreover, the semi-simple and nilpotent part of the Chevalley-Jordan decomposition of $\widetilde{U}$ are also derivations of $\mu$. Let us denote by $ \widetilde{U}_s$ the semi-simple part of $\widetilde{U}$. Since $\K$ is algebraically closed field, $ \widetilde{U}_s$
 is a diagonalizable endomorphism. Let $\lambda_1,\cdots,\lambda_n$ the set of eigenvalues. If we consider the Jordan-Chevalley decomposition of $U$, $U=U_s+U_n$, then $\widetilde{U_s}=\widetilde{U}_s$. Assume now that $\mu$ is rigid. We shall show in a first step that $U_s$ and $U_n$ belong to $ \mathcal{A}d_\mu$. Assume that $U_s \notin  \mathcal{A}d_\mu$.  Then $U_s$ is a non inner semisimple derivation of $\mu$. There exists $X_0 \in V$ such that $U=ad X_0$. Let $\{X_0,X_1,\cdots,X_{n-1}\}$ be a basis of eigenvectors of $U_s$ and $\{0,\lambda_1, \cdots,\lambda_{n-1}\}$ the corresponding eigenvalues. By hypothesis, we can also assume that these eigenvalues are non negative. Let us consider the deformation $\mu_\varepsilon$ of $\mu$ defined by 
 $$\mu_\varepsilon(X_0,X_i)=\mu(X_0,X_i)+\varepsilon U_s(X_i)$$
 and
 $$\mu_\varepsilon(X_i,X_j)=\mu(X_i,X_j)$$
for $1 \leq i < j \leq n-1.$ 
This define a non trivial deformation of $\mu$. Since $\g$ is rigid, we have a contradiction. Then $U_s \in \mathcal{A}d_\mu$ and $\mathcal{A}d_\mu$ is a split Lie algebra that is contains the semisimple and the nilpotent part of its elements. Now we can prove that $\mathcal{A}d_\mu$ is algebraic. Let us consider a semisimple element $U$ of $\mathcal{A}d_\mu$. If $\{\lambda_0=0,\lambda_1, \cdots,\lambda_{n-1}\}$  is the set of eigenvalues, its replica $\widetilde{U}$ is a semisimple element of $Der(\mu)$ whose eigenvalues  $\{\rho_0,\rho_1, \cdots,\rho_{n-1}\}$ satisfy $p_0\rho_0+p_1\rho_1+ \cdots+p_{n-1}\rho_{n-1}=0$ with
$(p_0,p_1,\cdots,p_{n-1}) \in \Lambda$ where $\Lambda$ is the subset of $\Z^n$ whose elements  satisfy
$p_0\lambda_0+ p_1\lambda_1 + \cdots+ p_{n-1}\lambda_{n-1} = 0.$ As abode, with such derivation we define an infinitesimal deformation of $\mu$. Since $\mu$ is rigid, the replica is in $\mathcal{A}d_\mu$ and this linear Lie algebra is algebraic.

 \section{Structure theorems of rigid Lie algebras}
 Assume that $\g$ is solvable with a trivial center. If $\g$ is rigid, then $ \mathcal{A}d_\mu$ is algebraic and we have the decomposition
$$ \mathcal{A}d_\mu= \n \oplus \frak{a}$$
where all the elements of $\n$ are nilpotent and $\frak{a}$ is an abelian subalgebra whose elements are semisimple. 
There exists $X_1,\cdots,X_{n-r},T_1,\cdots,T_r \in V$ such that $\{adX_1,\cdots,adX_{n-r}\}$ is a basis of $\n$ and thus $adX_i$ is nilpotent for $i=1, \cdots,n-r$ and $\{adT_1,\cdots,adT_r\}$ is a basis of $\frak{a}$ and thus $adT_i$ is semisimple for $i=1,\cdots,r$. Since $\frak{a}$ is abelian, we can assume that the endomorphisms $adT_i$ are diagonal. Let us denote by $\g_{\n}$ the subalgebra of $\g$ generated by $\{X_1,\cdots,X_{n-r}\}$ and $\g_{\frak{a}}$ the subalgebra generated by $\{T_1,\cdots,T_r \}$. Then we have

\begin{proposition}
 Let $\g$ be a finite dimensional solvable rigid $\K$ Lie algebra with a trivial center. Then $\g$ admits the decomposition
$$\g=\g_{\n} \oplus \g_{\frak{a}}$$
where $\g_{\n} $ is the nilradical of $\g$ and $\g_{\frak{a}}$ a maximal abelian subalgebra of $\g$ whose elements $T$ are such that $ad T$ is semisimple.
 \end{proposition}
\pf Since the operators $adX_i$ are nilpotent,  $\g_{\n} $ is a nilpotent subalgebra of $\g$. From the decomposition of $ \mathcal{A}d_\mu$, it is the maximal nilpotent ideal of $\g$. Let us note also that, sometimes, $\g_{\frak{a}}$  is called a maximal torus of $\g$. This may lead to some confusion since $\g$  is not necessarily algebraic. However, $\frak{a}$  is a maxiamal algebraic torus of the algebraic Lie algebra $ \mathcal{A}d_\mu$.

\medskip

\noindent{\bf Example.}
Let us consider the $2$-dimensional rigid Lie algebra
$$\g_2=
\left\{
\begin{pmatrix}
  x	&   x &y \\
   0   &  x & 0 \\
   0 & 0 & 0
\end{pmatrix}, \ x, y \in \C
\right\}.
$$
If 
$$X=
\begin{pmatrix}
  1	&   1 &0 \\
   0   &  1 & 0 \\
   0 & 0 & 0
\end{pmatrix}
, \ \ Y=
\begin{pmatrix}
  0	&   0 &1 \\
   0   &  0 & 0 \\
   0 & 0 & 0
\end{pmatrix}
$$
then $[X,Y]=Y$ and $ \mathcal{A}d_\mu=\{U,V\}$ with $U=adX,V=ad Y$. We have the decomposition
$$ \mathcal{A}d_\mu= \n \oplus \frak{a}$$
with $\n=\C\{V\}$ and $\frak{a}=\C\{U\}$. The corresponding decomposition of $\g$ is
$$\g=\C\{Y\} \oplus \C\{X\}.$$
We can see, in this example, that the elements of $\g_{\frak{a}}$ are not, in this case, semisimple.

\medskip

From this decomposition of $ \mathcal{A}d_\mu$ or of $\g$, Ancochea and Goze established an interesting criterium of rigidity. Let $\g= \g_{\frak{a}} \oplus \g_{\frak{n}}$ be the  decomposition of the rigid Lie algebra $\g$. 
\begin{definition} A vector $T_0 \in \g_{\frak{a}}$ is called regular if $\dim \ker (ad T_0) \leq \dim \ker (ad T)$ for any $T \in \g_{\frak{a}}$.
\end{definition}
 Moreover, since the elements of $\mathfrak{a}$ are semisimple and commuting, there exists a basis $\{T=X_1,X_2,\cdots,X_n\}$ of 
$\mathfrak{g}$ of eigenvectors for all the diagonal endomorphisms $ad T \in \frak{a}$, $T \in  \g_{\frak{a}}$. Let be $T \in \g_{\frak{a}}$. Then $$[T,X_i]=\lambda_i(T)X_i,$$
and the linear function $\lambda_i \in   \g_{\frak{a}}^*$ satisfy the relationswe
$$\lambda_i(T)+\lambda_j(T)=\lambda_k(T)$$ as soon as $[X_i,X_j]$ is an eigenvector corresponding to $\lambda_k(T)$. 
Let us denote by $S(T)$ the linear system whose equations are $$(x_i+x_j-x_k)=0$$  when  $C_{i,j}^k\neq 0.$
In particular, the linear system associated with the roots
$$\lambda_i(T)+\lambda_j(T)=\lambda_k(T)$$
is a subsystem of $S(T)$.

\begin{theorem}\cite{A.G-rank} Let $\g$ be a rigid solvable Lie algebra whose center is trivial. Then for any regular vector $T_0 \in \mathfrak{t}$, one has
$${\rm rank}(S(T_0))=\dim \mathfrak{n}-1.$$
\end{theorem}
An important example of such algebras are the Borel subalgebra of a semi-simple Lie algebras. Let us note also that this theorem permits to constructuct rigid Lie algebras without cohomological criterium. See, for example, \cite{ACO}.

 From the rank theorem, for any $X \in \g_{\n}-[\g_\n,\g_n]$, then $[X, \g_\mathfrak{a}] \neq 0.$ We deduce that $\g_\n= [\g_\n,\g_\n]$. We obtain
 
 \begin{proposition}
 Let $\g= \g_\mathfrak{a} \oplus \g_\n$ be a solvable rigid Lie algebra with a trivial center. Then the nilradical $\g_\n$ is the nilradical of $\g$ and it is an algebraic nilpotent Lie algebra.
 \end{proposition}
 \pf In fact, for any Lie algebra $\g$, its derived subalgebra $[\g,\g]$ is algebraic.

 For example, the "Heisenberg" Lie algebra  
 $$\frak{n}_1=\left\{ 
\left(
\begin{array}{llll}
x_1+x_2 & x_1+x_2 & 0 & x_1 \\
x_1+x_2 & x_1+x_2 & 0 & x_2 \\
 x_1& x_1+2 x_2 & 0 & x_3 \\
0 & 0 & 0 & 0 
\end{array}
\right), \ x_1,x_2,x_3 \in \C \right\} $$
cannot be the nilradical of a rigid Lie algebra.  On the other hand, the (good) Heisenberg Lie algebra
$$
\frak{n_2}= 
\left\{\begin{pmatrix}
   0   & x_1 & x_3   \\
    0  &  0 & x_2 \\
    0 & 0 & 0
\end{pmatrix}, \ x_1,x_2,x_3 \in \C \right\}$$
is the nilradical of a $5$-dimensional rigid Lie algebra whose multiplication is given by
$$
\left\{
\begin{array}{ll}
   \lbrack T_1, X_i \rbrack =X_i, i=1,3,    &   \lbrack T_2, X_i \rbrack=X_i, i=2,3, \\
  \lbrack X_1, X_2 \rbrack =X_3  &   
\end{array}
\right.
$$

\medskip

\noindent{\bf Remark.} Let $\g= \g_\mathfrak{a} \oplus \g_\mathfrak{n}$ the decomposition of a solvable rigid Lie algebra. We call root of $\g$, a non zero linear form $\alpha \in \g_\mathfrak{a}^*$ such that the linear space
$$\g_\alpha=\{ X \in \g_\n, [T,X]=\alpha (T) X\}$$
is not $\{0\}$. If we denote by $\Delta$ the set of roots, we have the decomposition
$$\g= \g_\mathfrak{a} \oplus_{\alpha \in \Delta} \g_\alpha.$$
When $\g$ is the Borel subalgebra of a semisimple Lie algebra, then $\Delta$ is the set of positive roots. As in this example, we can consider the subset of positive roots $\Pi$ of $\Delta$ and the nilradical $\g_\n$ admit a $\Delta$-grading. Recall that all the gradings of filiform Lie algebras are described in \cite{BGR}.

\noindent{\bf Examples.}

1. Let us consider the $5$-dimensional Lie algebra given by
$$
\left\{
\begin{array}{ll}
   \lbrack T_1, X_i \rbrack =X_i, i=1,3,    &   \lbrack T_2, X_i \rbrack=X_i, i=2,3, \\
  \lbrack X_1, X_2 \rbrack =X_3  &   
\end{array}
\right.
$$
The vector $T=T_1+T_2$ is regular and $S(T)$ is the linear system
$$
\left\{
\begin{array}{l}
t_2 +x_i=x_i, \ i=2,3, \\
x_1+x_2=x_3.  
\end{array}
\right.
$$
We have ${\rm rank}(S(T))=2=\dim \n-1$. Here $\n$ is the $3$-dimensional Heisenberg algebra and $\g$ is a Borel subalgebra of $sl(3)$.

2. Let us consider the $8$-dimensional Lie algebra given by
$$
\left\{
\begin{array}{ll}
   \lbrack T_1, X_i\rbrack =X_i, i=1,2,3,4,    &   \lbrack T_1, X_5 \rbrack = 2X_5 \\
     \lbrack T_2, X_i\rbrack =X_i, i=2,3,5,    &   \lbrack T_3, X_3 \rbrack = X_3 \\
  \lbrack T_3, X_4 \rbrack =-X_4,    & \lbrack X_1, X_2 \rbrack = \lbrack X_3, X_4 \rbrack =X_5    
\end{array}
\right.
$$
The vector $T_1$ is regular and $S(T)$ is the linear system
$$
\left\{
\begin{array}{l}
t_2 +x_i=x_i, \ i=2,4,5, \\
t_3 +x_i=x_i, \ i=3,4, \\
x_1+x_2=x_5, \\
x_3+x_4=x_5.  
\end{array}
\right.
$$
We have ${\rm rank}(S(T))=4=\dim \n-1$. But $\g$ is not a Borel Lie algebra.

\end{document}